\documentclass{amsart}
\usepackage{graphicx, color}
\newtheorem{theorem}{Theorem}[section]
\newtheorem{lem}[theorem]{Lemma}
\theoremstyle{remark}

\numberwithin{equation}{section}

\newcommand{\Real}{\mathbb R}

\begin{document}

\title[Lack of null controllability for the heat equation]{An elementary proof of the lack of null controllability for the heat equation on the half line}
\author{Konstantinos Kalimeris \and T\"urker \"Ozsar{\i}}
\address{Department of Applied Mathematics and Theoretical Physics, University of Cambridge, Cambridge, UK}%
\email{kk364@cam.ac.uk}
\address{Department of Mathematics, Izmir Institute of Technology, Izmir, Turkey}%
\email{turkerozsari@iyte.edu.tr}%
\maketitle

\begin{abstract}
  In this note, we give an elementary proof of the lack of null controllability for the heat equation on the half line by employing the machinery inherited by the unified transform, known also as the Fokas method. This approach also extends in a uniform way to higher dimensions and different initial-boundary value problems governed by the heat equation, suggesting a novel methodology for studying problems related to controllability.
\end{abstract}

\section{Introduction}
The Uniform Transform Method (UTM), also known as the Fokas method, is a powerful tool for obtaining solutions  of initial - (inhomogeneous) boundary-value problems.  This method was first introduced in \cite{Fok97} for the analysis of initial-boundary value problems for integrable nonlinear partial differential equations (PDEs).  However, in later works it was proven to produce novel results for a general class of linear PDEs; see \cite{Fok02,Fok08}. Recently researchers utilized the UTM to produce rigorous wellposedness results in Sobolev and Bourgain spaces for dispersive PDEs; see for instance \cite{Him17} and \cite{Ozs19} for the local and global wellposedness analysis of nonlinear Schr\"{o}dinger type PDEs and \cite{Him19} for a similar analysis on the Korteweg-de Vries equation.

To date, there is no work on the boundary controllability of PDEs that utilizes the advantages of the UTM. This method has two basic elements: (i) the so-called \emph{global relation}, an identity that relates the initial datum and a suitable time transform of known and unknown boundary values, and (ii) the \emph{integral representation} of the solution.  We illustrate a new methodology by making use of these two elements in order to provide an elementary proof of the lack of null controllability for the heat equation on the half line.

To this end, let us consider the following canonical initial-boundary value problem:
\begin{eqnarray}
  u_t &=& u_{xx},\quad x\in \mathbb{R}_+, \ \ t\in (0,T), \label{maineq1}\\
  u(x,0) &=& u_0(x),\quad x\in \mathbb{R}_+, \label{maineq2}\\
  u(0,t) &=& g(t),\quad t\in (0,T). \label{maineq3}
\end{eqnarray}  We say \eqref{maineq1}-\eqref{maineq3} is \emph{null controllable} in $[0,T]$ if given $u_0\in L^2(\mathbb{R}_+)$ there is $g\in L^2(0,T)$ such that $u(x,T)\equiv 0$.

It is well known that the above property does not hold for \eqref{maineq1}-\eqref{maineq3} for those solutions in $C([0,T];L^2(\mathbb{R}_+))$; see for example \cite{MicZua00} for a proof of this result. Our goal is to provide an alternate, yet very short proof of this fact.  More precisely, we  prove the following theorem.
\begin{theorem}\label{mainthm}
  There exists $u_0\in L^2(\mathbb{R}_+)$ such that $u(x,T)\not\equiv 0$ for any $g\in L^2(0,T)$ if $u\in C([0,T];L^2(\mathbb{R}_+))$ and it solves \eqref{maineq1}-\eqref{maineq3}.
\end{theorem}

\subsection*{Orientation} In Section 2, we provide a proof of Theorem \ref{mainthm} via the global relation. In Section 3, we extend Theorem \ref{mainthm} to the $N$-dimensional half space by outlining the straightforward and simple extension of the proof presented in Section 2 to $N$ dimensions. In Section 4, we discuss alternative pathways through the Fokas method, introducing also a characterisation for the null-controllability problem on the finite interval. In Section 5, we discuss the main results of this work, as well as its future implications.

\section{Proof of Theorem \ref{mainthm}}

By introducing the half-line Fourier $x$-transform, namely
\begin{equation} \label{half-F-1}
\hat{f}(\lambda)=\int_0^{\infty}e^{-i\lambda x}f(x)dx,  \qquad \text{Im}\lambda\leq 0,
\end{equation}
and
\begin{equation*} 
\hat{F}(\lambda,t )=\int_{0}^{\infty }e^{-i\lambda x}F(x,t)dx,  \hphantom{2a}  \text{Im}\lambda\leq 0 ,   
\end{equation*}
as well as  the $t$-transform
\begin{equation} \label{t-tr-1}
\tilde f(\lambda,t)=\int_0^{t}e^{\lambda \tau}f(\tau)d\tau, \ \ t>0, \  \lambda \in \mathbb{C},
\end{equation}
 the global relation for \eqref{maineq1}-\eqref{maineq3}, given by the Fokas method (equation (12) in \cite{Fok08}) can be written in the following form:
\begin{equation}  \label{GR-1-hl}
e^{\lambda^2t}\hat{u}(\lambda,t)=\hat{u}_0(\lambda)-\tilde{r}(\lambda^2,t)-i\lambda\tilde{g}(\lambda^2,t),\qquad \text{Im}\lambda\leq 0,
\end{equation}
where $ r(t)=u_x(0,t)$ and $ g(t)=u(0,t),\ \  t>0 $.  For matters of completeness we derive here the global relation using the half-Fourier transform. Indeed, through integration by parts we obtain
\begin{align*}
\hat{u}_t(\lambda ,t)&=\int_0^{\infty}e^{-i\lambda x}u_t(x,t) dx=\int_0^{\infty}e^{-i\lambda x}u_{xx}(x,t)dx\\
&=u_x(x,t)e^{-i\lambda x}\big|_{x=0}^{\infty}+i\lambda u(x,t)e^{-i\lambda x}\big|_{x=0}^{\infty}-\lambda^2 \hat{u}(\lambda,t).
\end{align*}
Thus,
$$\hat{u}_t+\lambda^2\hat{u}=-r(t)-i\lambda g(t).$$
Integrating the above ordinary differential equation we obtain
$$\hat{u}e^{\lambda^2 t}=\hat{u}_0-\int_0^t e^{\lambda^2\tau} [r(\tau )+i\lambda g(\tau )]d\tau,$$
which is \eqref{GR-1-hl}.

 Applying the condition $u(x,T)\equiv 0$ in \eqref{GR-1-hl}, we obtain that
\begin{equation}  \label{GR-2-hl}
0=\hat{u}_0(\lambda)-\tilde{r}(\lambda^2,T)-i\lambda\tilde{g}(\lambda^2,T),\qquad \text{Im}\lambda\leq 0.
\end{equation}
Letting $\lambda\to -\lambda$ in \eqref{GR-2-hl} and subtracting the resultant expression (which is valid for $\text{Im}\lambda\ge 0$) from \eqref{GR-2-hl} we obtain the following equation:
\begin{equation}\label{main-eq-hl}
2i\lambda\tilde{g}(\lambda ^{2},T)=\hat{u}_{0}(\lambda )-\hat{u}_{0}(-\lambda ), \qquad \lambda\in\Real.
\end{equation}

 Let $0\not\equiv u_0\in L^1\cap L^2(\mathbb{R}_+)$. 
Employing this assumption  in \eqref{main-eq-hl} along with the definition of $\tilde{g}$, we obtain the following uniform bound for some $M>0$:
\begin{equation}\label{impineq}
  {\left|\int_0^Te^{\lambda ^2t}g(t)dt\right|}=\left|\frac{1}{2\lambda}[\hat{u}_{0}(\lambda )-\hat{u}_{0}(-\lambda )]\right|{<M}, \qquad \lambda^2>1.
\end{equation}
Then ${g\equiv 0}$ due to the Lemma \ref{Yoslem}, below.

It is clear that if $g\equiv 0$, then ${\hat{u}_{0}(\lambda )}={\hat{u}_{0}(-\lambda )}$ for all $\lambda\in \mathbb{R}$, which would contradict with the assumption that $0\not\equiv u(0)=u_0$.

\begin{lem}\label{Yoslem}(\cite{Yosida}, page 167, Lemma 2) Let $g\in L^2(0,T)$. If there is $M>0$ such that $\left|\int_0^Te^{\alpha t}g(t)dt\right|<M$ for every $\alpha>1$, then $g\equiv 0$.
\end{lem}
 We note that the proof in \cite{Yosida} is given for $g$ being a continuous function;  the proof  extends to $L^2$ functions via density, namely $g$ is vanishing almost everywhere.
 
\section{The $N$-dimensional half space}

In this section we extend Theorem \ref{mainthm} to the higher dimensional half space $\mathbb{R}^N_+=\mathbb{R}^{N-1}\times  \mathbb{R}_+$, $\ N>1$ (see also \cite{MicZua01}). The methodology we used previously for the proof of Theorem \ref{mainthm} provides a straightforward path  to  study the (lack of) null controllability for
  \begin{eqnarray}
  u_t &=& \Delta u,\quad x=(x',x_N)\in \mathbb{R}^N_+, \ \  t\in (0,T), \label{hmaineq1}\\
  u(x,0) &=& u_0(x),\quad x\in \mathbb{R}^N_+, \label{hmaineq2}\\
  u(x',0,t) &=& g(x',t),\quad x'\in \mathbb{R}^{N-1}, \ \ t\in (0,T). \label{hmaineq3}
\end{eqnarray}


The relevant result can be obtained by using half space Fourier $x$-transform $$\displaystyle\hat{u}(\lambda)\doteq\int_{\mathbb{R}^{N-1}}\int_0^\infty e^{-i\lambda\cdot x}u(x)dx_ndx', \qquad \lambda=(\lambda',\lambda_N)\in\mathbb{R}^{N-1}\times \mathbb{C},\quad \text{Im}\lambda_N\le 0$$ and applying Fokas's method only to the last variable $x_N$.  Indeed, half space Fourier transform yields the global relation
\begin{equation}\label{HSGR}
  e^{|\lambda|^2t}\hat{u}(\lambda,t)=\hat{u}_0(\lambda)-\tilde{h}(\lambda,t)-i\lambda_N\tilde{g}(\lambda,t), \qquad  \text{Im}\lambda_N\le 0,
\end{equation} where
\begin{equation}\label{boundaryvalues}
  \tilde{g}(\lambda,t)\doteq\int_0^te^{|\lambda|^2s}\widehat{g^{x'}}(\lambda',s)ds \ \ \ \text{ and } \ \ \  \tilde{h}(\lambda,t)\doteq\int_0^te^{|\lambda|^2s}\widehat{h^{x'}}(\lambda',s)ds,
\end{equation} with $h(x',t)\doteq u_{x_N}(x',0,t)$ and $\widehat{g^{x'}}$, $\widehat{h^{x'}}$  denoting Fourier transforms of $g$ and $h$ with respect to $x'$.


The proof of the lack of null controllability for solutions in the class $C([0,T];L^2(\mathbb{R}^N_+))$ follows the exact same steps with the proof of Theorem \ref{mainthm}. Hence, \eqref{impineq} is now replaced with
\begin{equation}\label{impineq2}
  {\left|\int_0^Te^{\lambda_N^2t}F(\lambda',t)dt\right|}=\left|\frac{1}{2\lambda_N}[\hat{u}_{0}(\lambda',\lambda_N )-\hat{u}_{0}(\lambda',-\lambda_N)]\right|{<M}, \ \ \ \ \lambda_N^2>1,
\end{equation} where $F(\lambda',t):=e^{|\lambda'|^2t}\widehat{g^{x'}}(\lambda',t)$.  Applying Lemma \ref{Yoslem} for each fixed $\lambda'\in \mathbb{R}^{N-1}$, we conclude that $F\equiv 0$, which in turn implies that $g\equiv 0$.

\section{Alternative Pathways}

In this section, we provide an alternative pathway to obtain a proof of Theorem \ref{mainthm} via the integral representation of the Fokas method. Furthermore, this pathway provides a characterisation of the control for the finite interval problem given in \eqref{main-eq-fi}. In this sense it suggests a more general viewpoint on studying controllability problems through this methodology.

\subsection*{The Half Line}
The integral representation of the  solution of \eqref{maineq1}-\eqref{maineq3} given by the Fokas method (equation (16) in \cite{Fok08}) takes the form:
  \begin{align}\label{repform}
u(x,t)&=\frac{1}{2\pi}\int_{-\infty }^{\infty }e^{i\lambda x-\lambda ^{2}t}\hat{u}_{0}(\lambda )d\lambda \notag \\
&-\frac{1}{2\pi}\int_{\partial D^{+}}e^{i\lambda x-\lambda ^{2}t}
\left[ 2i\lambda \tilde{g}(\lambda ^{2},t)+\hat{u}_{0}(-\lambda )    \right] d\lambda,
\end{align}
 where $\partial D^+$ is depicted in Figure 1.
 \begin{figure}[!ht]
\label{fig}
\centering
\includegraphics[width=0.4\linewidth]{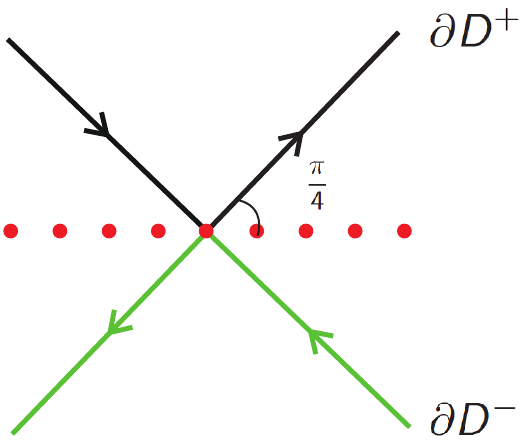}
\caption{The contours $\partial D^\pm$.}
\end{figure}

By applying $u(x,T)\equiv 0$, deforming $\partial D^+$ to the real line and taking the inverse Fourier transform of both sides in the resultant expression, we obtain \eqref{main-eq-hl}. Then, the proof of Theorem \ref{mainthm} follows by the exact same arguments of Section 2.

\subsection*{The Finite Interval} It is well known that the null controllability is true, for instance in $C([0,T];L^2(\Omega))$, if one replaces the infinite domain $\mathbb{R}_+$ by the finite one $(0,L)$.  Here, we wish to give a characterization of the set of suitable boundary controllers, say acting at the right Dirichlet boundary condition, using the integral representation obtained from the Fokas method.
Thus, we consider the following problem:
\begin{eqnarray}\label{heat-fi}
 \begin{cases}
  u_t =  u_{xx}, \qquad &x\in (0,L), \ t\in (0,T),\\
  u(0,t)=0, \quad u(L,t)=h(t), \qquad  &t\in (0,T)\\
  u(x,0)=u_0(x), \qquad &x\in (0,L)
 \end{cases}
\end{eqnarray} and the goal is to find a sufficient condition for the boundary controller $h$ so that it steers the given initial datum $u_0$ to $u_T\equiv 0$ at $t=T$.

In analogy with the half line problem, one introduces the following Fourier $x$-transform where the integral is taken over the given spatial domain $(0,L)$:
\begin{equation}  \label{fourthree}
\hat{u}(\lambda ,t)=\int_{0}^{L}e^{-i\lambda x}u(x,t)dx, \qquad \lambda \in\mathbb{C}.
\end{equation}  Then, the corresponding global relation (equation (2.10) in \cite{Fok08}) for the above problem evaluated at $t=T$ becomes
\begin{equation}\label{globalrel}
  0=\hat{u}_0(\lambda)+i\lambda e^{-i\lambda L}\tilde{h}(\lambda^2,T)-\tilde{g}_1(\lambda^2,T)+e^{-i\lambda L}\tilde{h}_1(\lambda^2,T), \qquad \lambda\in \mathbb{C},
\end{equation} 
with $g_1(t)=u_x(0,t)$, $h_1(t)=u_x(L,t)$, and ${h}(t)=u(L,t)$. 

Similarly, the integral representation of the solution (equation (2.6) in \cite{Fok08}) evaluated at $t=T$ becomes
\begin{multline}\label{int-rep-fi}
0=u(x,T)=\frac{1}{2\pi}\int_{-\infty }^{\infty }e^{i\lambda x-\lambda ^{2}T}\hat{u}_{0}(\lambda )d\lambda 
-\frac{1}{2\pi}\int_{\partial D^{+}}e^{i\lambda x-\lambda ^{2}T}
 \tilde{g}_{1} (\lambda^2,T)  d\lambda   \\  
-\frac{1}{2\pi}\int_{\partial D^{-}}e^{-i\lambda (L-x)-\lambda ^{2}T}
\left[ \tilde{h}_{1}(\lambda^2,T)+i\lambda \tilde{h}(\lambda^2,T)   \right] d\lambda , 
\end{multline}
for all $x\in(0,L)$, where the contours $\partial D^\pm$ are depicted in Figure 1.

We next utilise the standard approach of Fokas method:
Using the invariances of the global relation under the transformation $\lambda\mapsto-\lambda$, the unknown boundary transforms ($\tilde{g}_1$ and $\tilde{h}_1$) can be eliminated from the integral representation (see equation (32) in \cite{Fok08}). Through short and straightforward calculations, and by employing the definition of $\tilde{h}$, equation \eqref{int-rep-fi} yields the following relation:
\begin{equation}\label{main-eq-fi}
 \int_{\partial D_+}R(\lambda;x,T,L)d\lambda+\int_{\partial D_-}R(\lambda; x,T,L)d\lambda=U_0(x;T),  \ \ \  \forall  \ x\in (0,L),
  \end{equation}
where the integrand $R(\lambda;x,T,L)$ is given by
\begin{equation}\label{R-fi}
 R(\lambda;x,T,L):=\frac{i}{\pi}\frac{\lambda e^{i\lambda x-\lambda^2T}}{e^{i\lambda L}-e^{-i\lambda L}}\left[\int_0^Te^{\lambda^2 s}h(s)ds\right]
\end{equation}
and the known $U_0(x;T)$ is given by
\begin{multline}\label{U-fi}
U_0(x;T)=  \frac{1}{2\pi}\int_{-\infty}^\infty e^{i\lambda x-\lambda^2T}\hat{u}_0(\lambda)d\lambda 
\\ -\frac{1}{2\pi}\int_{\partial D_+}e^{i\lambda x-\lambda^2T}\left[\frac{e^{i\lambda L}\hat{u}_0(\lambda)-e^{-i\lambda L}\hat{u}_0(-\lambda)}{e^{i\lambda L}-e^{-i\lambda L}}\right]d\lambda
  \\-\frac{1}{2\pi}\int_{\partial D_-}e^{-i\lambda (L-x)-\lambda^2T}\left[\frac{\hat{u}_0(\lambda)-\hat{u}_0(-\lambda)}{e^{i\lambda L}-e^{-i\lambda L}}\right]d\lambda,
\end{multline}
with the contours $\partial D^\pm$  depicted in Figure 1, and the red dots denoting the zeros of $\exp (i\lambda L)-\exp (-i\lambda L)$ on the real axis.


Thus, we obtain the following characterization for the problem of null controllability:
The problem \eqref{heat-fi} is null controllable at time $t=T$ if and only if there exists $h=h(t)$ which satisfies \eqref{main-eq-fi}.

  \section{Discussion}
  
  In this work we analyse a family of null-controllability problems governed by the heat equation, using the machinery provided by the Fokas method. In this connection we make the following three remarks:
  \begin{itemize}
 \item It is straightforward but more technical to generalise the proof of Theorem \ref{mainthm}, so that one constructs a function $u_0$ satisfying Theorem \ref{mainthm}, with $u_0\in L^2(\mathbb{R}_+)$, but not necessarily $u_0\in L^1(\mathbb{R}_+)$.
 \item  The methodology appearing in the current work can be applied to boundary value problems of higher dimensions such as $(\mathbb{R}_+)^N, \ N >1$, where all the spatial coordinates are positive. The relevant proof, which will be presented elsewhere, is based on the analysis of the Fokas method presented in \cite{Fok02} for the case of $N=2$, namely the quarter plane.
\item  If $u_0\in L^2(\mathbb{R}_+)$ and $g\in L^2(0,T)$, then \eqref{maineq1}-\eqref{maineq3} possesses a solution $u\in C([0,T];L^2(\mathbb{R}_+))$ in the transposition sense, and moreover this solution can be represented as in \eqref{repform}.  Therefore, Theorem \ref{mainthm} concerns such solutions.  If the condition $u\in C([0,T];L^2(\mathbb{R}_+))$ is removed, then one can recover the null controllability in a larger class of solutions.  This was proved in \cite{Ros2000} for the linearized KdV, heat, and Schr\"odinger equations.
\end{itemize}

The Fokas method provides the basic tools which are needed for the extension of the methodology introduced in the current work to linear PDEs, other than the heat equation. Indeed, one could obtain the Global Relation of the initial and boundary conditions, as well as the Integral Representation of the solution for problems which are posed on the half line and the finite interval and satisfy evolution equations where the rhs of \eqref{maineq1} is substituted by a higher order linear differential operator with constant coefficients (see \cite{Fok08}).
The possibility of applying this methodology to null-controllability problems governed by other linear evolution PDEs  is currently under investigation.

\section*{Acknowledgment}
The authors wish to thank A.S. Fokas (University of Cambridge), whose prolific works are an endless source of inspiration. KK acknowledges funding by EPSRC. T\"O research is funded by TUBITAK 1001 Grant {\#}117F449.

\bibliographystyle{amsplain}

\begin{thebibliography}{10}

\bibitem{Fok97}
Athanassios~S. Fokas, \emph{A unified transform method for solving linear and
  certain nonlinear pdes}, Proceedings of the Royal Society of London. Series
  A: Mathematical, Physical and Engineering Sciences \textbf{453} (1997),
  no.~1962, 1411--1443.

\bibitem{Fok02}
\bysame, \emph{A new transform method for evolution partial differential
  equations}, IMA Journal of Applied Mathematics \textbf{67} (2002), no.~6,
  559--590.

\bibitem{Fok08}
\bysame, \emph{A unified approach to boundary value problems}, CBMS-NSF
  Regional Conference Series in Applied Mathematics, vol.~78, Society for
  Industrial and Applied Mathematics (SIAM), Philadelphia, PA, 2008.
  \MR{2451953}

\bibitem{Him17}
Athanassios~S. Fokas, A.~Alexandrou Himonas, and Dionyssios Mantzavinos,
  \emph{The nonlinear {S}chr\"{o}dinger equation on the half-line}, Trans.
  Amer. Math. Soc. \textbf{369} (2017), no.~1, 681--709. \MR{3557790}

\bibitem{Him19}
A.~Alexandrou Himonas, Dionyssios Mantzavinos, and Fangchi Yan, \emph{The
  {K}orteweg--de {V}ries equation on an interval}, J. Math. Phys. \textbf{60}
  (2019), no.~5, 051507, 26. \MR{3947621}

\bibitem{MicZua00}
Sorin Micu and Enrique Zuazua, \emph{On the lack of null-controllability of the
  heat equation on the half-line}, Trans. Amer. Math. Soc. \textbf{353} (2001),
  no.~4, 1635--1659. \MR{1806726}

\bibitem{MicZua01}
\bysame, \emph{On the lack of null-controllability of the heat equation on the
  half space}, Portugaliae Mathematica \textbf{58} (2001), no.~1, 1--24.

\bibitem{Ozs19}
T{\"u}rker {\"O}zsar{\i} and Nermin Yolcu, \emph{The initial-boundary value
  problem for the biharmonic {S}chr{\"o}dinger equation on the half-line},
  Commun. Pure Appl. Anal. \textbf{18} (2019), no.~6, 3285--3316.

\bibitem{Ros2000}
Lionel Rosier, \emph{Exact boundary controllability for the linear
  {K}orteweg-de {V}ries equation on the half-line}, SIAM J. Control Optim.
  \textbf{39} (2000), no.~2, 331--351. \MR{1788062}

\bibitem{Yosida}
K\^{o}saku Yosida, \emph{Functional analysis}, sixth ed., Grundlehren der
  Mathematischen Wissenschaften [Fundamental Principles of Mathematical
  Sciences], vol. 123, Springer-Verlag, Berlin-New York, 1980. \MR{617913}

\end{thebibliography}

%
%

\end{document}